\newtheorem{lemma}{Lemma}[section]
\newtheorem{prop}[lemma]{Proposition}
\newtheorem{theorem}[lemma]{Theorem}
\newtheorem{cor}[lemma]{Corollary}
\newtheorem{rem}[lemma]{Remark}
\newtheorem{defi}{Definition}

\newcommand{\lem}{\begin{lemma}}
\newcommand{\ma}{\end{lemma}}
\newcommand{\pro}{\begin{prop}}
\newcommand{\tion}{\end{prop}}
\newcommand{\theo}{\begin{theorem}}
\newcommand{\orem}{\end{theorem}}
\newcommand{\coro}{\begin{cor}}
\newcommand{\lar}{\end{cor}}
\newcommand{\remm}{\begin{rem}}
\newcommand{\ark}{\end{rem}}
\newcommand{\defin}{\begin{defi}}
\newcommand{\nion}{\end{defi}}

\newcommand{\kla}{\left ( }
\newcommand{\nach}{\rightarrow}
\newcommand{\mer}{\right ) }
\newcommand{\mitt}{\left | { \atop } \right.}
\newcommand{\for}{\begin{eqnarray*}}
\newcommand{\mel}{\end{eqnarray*}}

\newcommand{\nz}{{\rm  I\! N}}
\newcommand{\nen}{n \in \nz}
\newcommand{\ken}{k \in \nz}

\newcommand{\leb}{{\cal L}}

\newcommand{\ef}{{\cal F}}

\newcommand{\p}{\hspace{.05cm}}
\newcommand{\pl}{\hspace{.1cm}}
\newcommand{\pll}{\hspace{.3cm}}
\newcommand{\pla}{\hspace{1.5cm}}
\newcommand{\hz}{\vspace{0.5cm}}
\newcommand{\hoch}{ \vspace{-0.86cm}}

\newcommand{\Om}{\Omega}

\newcommand{\al}{\alpha}

\newcommand{\si}{\sigma}

\newcommand{\tet}{\theta}
\newcommand{\ttett}{\vartheta}

\newcommand{\eps}{\varepsilon}

\newcommand{\ds}{D_{\si}}

\newcommand{\lif}{\ell_{\infty}}

\newcommand{\noo}{\left \|}
\newcommand{\rrm}{\right \|}

\newcommand{\bet}{\left |}
\newcommand{\rag}{\right |}

\newcommand{\intt}{\int\limits}
\newcommand{\summ}{\sum\limits}
\newcommand{\limm}{\lim\limits}

\newcommand{\pq}{\pi_{pq}}
\newcommand{\pr}{\pi_{r1}}

\newcommand{\gc}{\tilde{C}}

\documentstyle[11pt]{article}
\begin{document}

\title{How many vectors are needed to compute (p,q)-summing
norms?\footnote{A preliminary
version of this paper occurred as "Absolutely summing norms with
n vectors"\cite{DJ2}}}

\author{Martin Defant \\ Marius Junge }

\date{}

\maketitle
%%%%%%%%%%%%%%%%%%%%%%%%%%%%
%%%%%%%%%%%%%%Introduction
%%%%%%%%%%%%%%%%%%%%%%%%%%%%%%

\newtheorem{llemma}{Lemma}
\newtheorem{pprop}[llemma]{Proposition}
\newtheorem{ttheorem}[llemma]{Theorem}

\begin{abstract}
\parindent0em
%%%%%%
\newcommand{\baselinestrech}{1.1}
\small\normalsize
\hspace{-1,9em} We show that for $q<p$ there exists an $\al < \infty$
such that
\[ \pq(T) \pl \le \pl c_{pq}\pl \pq^{[n^{\al}]}(T) \pla
\mbox{for all $T$ of rank $n$.}\]
Such a polynomial number is only possible if $q=2$ or $q<p$.
Furthermore, the growth rate is linear if $q=2$ or
$\frac{1}{q}-\frac{1}{p}>\frac{1}{2}$. Unless
$\frac{1}{q}-\frac{1}{p}=\frac{1}{2}$ this is also
a necessary condition. Based on similar ideas we prove that for $q>2$
the Rademacher cotype constant of a $n$-dimensional Banach space can be
determined with essentially $n (1+\ln n)^{c_q}$
many vectors.
\end{abstract}

\setcounter{section}{0}
\section*{Introduction}

In the local theory of Banach spaces the concept of summing operators
is of special interest. The presented paper is concerned with
the following problem raised up by T. Figiel.\hz

For given $1 \le q \le p \le \infty$ what is the best rate
$k_n$, such that
\vspace{2ex}

{\bf(*)\pla}\vspace{-5ex}
\for \pq(T)&\le& c \pl \pq^{k_n}(T) \mel
holds for all operators of rank $n$ and some constant $c$?\hz

In \cite{DJ} an observation of Figiel and Pelczynski
was generalized in showing
\[\pq(T) \pl \le \pl 3\pl \pq^{16^n}(T)  \]
for all q, p and all operators $T$ of rank $n$. This exponential growth
can not be improved in general. Figiel and Pelczynski also
showed that there is an operator $T: \lif^{2^n}\nach \ell_2^n$
(the Rademacher projection) such that for all $\ken$
\[ \pi_1^k(T) \pl \le \pl e\pl \sqrt{\frac{1+ \ln k}{n}}\pl
 \pi_1(T)\pl . \]
Recently, Johnson and Schechtman \cite{JOS} discovered that
for $p=q$ and $q \neq 2$ the rate can not be polynomial. More
precisely, every sequence satisfying (*) growth faster than
any polynomial, i.e.
\[ \limm_{n \nach \infty} k_n \p n^{-t} \pl=\pl \infty \]
for all $0<t<\infty$.
This phenomenon is related with the fact that
$L_p$ spaces don't have the polynomial approximation property,
which was proved by Bourgain.\hz

By far the nicest and most important result is Tomczak-Jaegermann's
inequality, namely
\vspace{2ex}

{\bf(1)\pla}\vspace{-5ex}
\for
\pi_2(T) &\le& \sqrt{2} \pl \pi_2^n(T)
\pla \mbox{for all $T$ of rank $n$.} \mel
In \cite{DJ} a certain type of
quotient formula was used to generalize Tomczak-Jaegermann's
inequality:
\[ \pi_{p2}(T) \pl \le \pl \sqrt{2} \pl \pi_{p2}^n(T)
\pla \mbox{for all $T$ of rank $n$.}\]
K\"onig and Tzafriri showed that
for all $2 < p< \infty$
\vspace{2ex}

{\bf(2)\pla}\vspace{-5ex}
\for
\pi_{p1}(T) &\le& c_p \pi_{p1}^n(T)
\pla \mbox{for all $T$ of rank $n$.} \mel
\hz

In contrast to the case $p=q$ we can show that for $q<p$ the
(p,q)-summing norm can be well-estimated by a polynomial number
of vectors.\hz

\begin{ttheorem}\label{polynom} Let $1 \le q \le p\p,r \le \infty$ with
$\frac{1}{q}=\frac{1}{p}+\frac{1}{r}$. Then for all operator T
of rank n one has

\[ \pq(T) \pl \le \pl c_r \left\{ \begin{array}{l@{\quad}l}
 \pq^n(T)              & for \pl 1 \le r < 2 \\[+0.2cm]
 \pq^{[n(1+\ln n)]}(T) & for \pl 2=r \\[+0.2cm]
 \pq^{[n^{r/2}]}(T )   & for \pl  2 < r<\infty \pl .
\end{array} \right.\]
\end{ttheorem}\hz

A very helpful tool in the proof of this theorem is again a quotient
formula for (p,q)-summing operators which allows a reduction to
the (probably worst) case $q=1$.\hz

\begin{ttheorem}\label{quot} Let $1\le q \le p \le \infty$ and
$1 \le r\le s\le q'$ with $\frac{1}{r}=\frac{1}{p}+\frac{1}{s}$.
Then for all operator $T: X \nach Y$ and $\nen$
\[ \pq^n(T) \pl = \pl \sup \{ \p \pr^n\kla TV\ds\mer \p \mid \p
V:\ell_{q'} \nach X,\pl \ds:\lif \nach \ell_{q'}, \pl \noo \si \rrm_s,
\noo V \rrm\p \le \p 1 \} \pl . \]
\end{ttheorem}\hz

For instance the first case of theorem 1 is a direct consequence
of theorem 2 and (2). In the other cases a crucial
observation of Jameson gives the link between limit orders and
number of vectors, see chapter 2. We are in debt to W.B. Johnson
for showing us Jameson's paper \cite{JAM}. There has always been a
quite
close connection between the theory of absolutely (p,q)-summing
operators
and the theory of cotype in Banach spaces. For instances, as a
consequence of
Tomczak-Jaegermann's inequality and it's generalization the gaussian
cotype
constant of an $n$-dimensional Banach space can be calculated with $n$
vectors.
This problem is still open in the case of Rademacher cotype. The
presented
technique allows us to reduce the number of vectors to the order $n (1+
\ln n)^{c_q}$
which indicates a positive solution for the Rademacher cotype.
Unfortunately, the constant $c_q$ tends to infinity as $q$ tends to
$2$.

\begin{ttheorem} \label{Rcot} Let $2<q<\infty$ and E a $n$ dimensional
Banach space.
For the Rademacher cotype constant one has
\[ C_q(Id_E) \pl \le \pl 2 \pl C_q^m (Id_E) \pl ,\]
where $m$ satisfies the following estimate for an absolute constant
$c_0$
\[ m \pl \le\pl n \pl (c_0\p (1+\ln n))^{\frac{1}{1-\frac{2}{q}}} \pl
.\]
\end{ttheorem}

Finally, we want to indicate that a linear growth is only possible if
$q=2$ or
$\frac{1}{q}-\frac{1}{p}\le\frac{1}{2}$.

\begin{ttheorem} \label{exam} Let $1 \le q \le p < \infty$, $q \le r
\le \infty$
with $\frac{1}{q}=\frac{1}{p}+\frac{1}{r}$. If $q \neq 2$ and
$2 < r \le \infty$. Then there exists $\al >1$ such that for all
sequences $k_n$ with
\[ \pi_{pq}(T) \pl \le \pl c \pl \pq^{k_n}(T) \pla
\mbox{for all T of rank n}\]
there exists a constant $\bar{c}$ with
\[ n^{\al} \pl \le \pl \bar{c} \pl k_n \pl . \]
\end{ttheorem}\hz

The constructed examples are very closely related to limit orders
of (p,q)-summing operators. In fact, it is well known that the identity
on $\ell_2^n$ yields an example for the proposition above as long as
$q>2$. In the case $q<2$ we intensively use the results of
Carl, Maurey and Puhl \cite{CMP} about Benett matrices.

%%%%$$$$$$$$$$$$$$$$$$$$$$$
%%%Preliminaries
%%%%%$$$$$$$$$$$$$$$$$$$$$$$

\section*{Preliminaries}
In what follows $c_0, c_1,$ .. always denote universal constants.
We use standard Banach space notation. In particular, the classical
spaces $\ell_q$ and $\ell_q^n$, $1\le q\le \infty$, $\nen$, are defined
in the usual way. By $\iota : \ell_q^n \nach \ell_p^n$ we denote the
canonical identity. Let $(e_k)_{\ken}$ be the sequence of unit vectors
in $\lif$.
For a sequence $\si = (\si_k)_{\ken} \in \lif$,
$\tau = (\tau_k)_{\ken} \in \lif$ we define
\[ \ds(\tau) \pl := \pl \summ_k \si_k \p \tau_k \p e_k \pl . \]
The standard reference on operator ideals is the monograph of Pietsch
\cite{PIE}. The ideals of linear bounded operators, finite rank
operators,
integral operators are denoted by $\leb$, $\ef$, ${\cal I}$. Here
the integral norm of $T \in {\cal I}(X,Y)$ is defined by
\[ \iota_1(T) \pl := \pl  \sup \{ \bet tr(ST)\rag \mid S \in \ef(Y,X),
 \pl \noo S \rrm \le 1\}\pl . \]
\hz

Let $1 \le q \le p \le \infty$ and $\nen$. For an operator $T \in
\leb(X,Y)$
the pq-summing norm of T with respect to $n$ vectors is defined by
\[ \pq^n(T) \pl := \pl
 \sup\left\{\p \kla \summ_1^n \noo Tx_k \rrm^p \mer^{1/p} \p \bet \pl
  \sup_{\noo x^* \rrm_{X^*}\le1} \kla \summ_1^n \bet\langle
  x_k,x^*\rangle \rag^q \mer^{1/q}
 \pl \le \pl 1\right.\p \right\} \pl .\]
An operator is said to be absolutely pq-summing, short pq-summing,
$(T \in \Pi_{pq}(X,Y))$ if
\[ \pq(T) \pl := \pl \sup_n \pq^n(T) \pl < \pl \infty \pl . \]
Then $(\Pi_{pq},\pq)$ is a maximal and injective Banach ideal (in the
sense
of Pietsch). As usual we abbreviate  $(\Pi_q,\pi_q) :=
(\Pi_{qq},\pi_{qq})$.
For further information about absolutely pq-summing operators we refer
to
the monograph of Tomczak-Jaegermann \cite{TOJ}. In particular, we would
like to mention an elementary observation of Kwapien, see \cite{TOJ}.
Let
$1\le q \le p \le \infty$, $1\le \bar{q} \le \bar{p} \le \infty$ with
$q \le \bar{q}$, $p \le \bar{p}$ and $\frac{1}{q}-\frac{1}{p}=
\frac{1}{\bar{q}}-\frac{1}{\bar{p}}$. Then one has for all $T$
\[ \pi_{\bar{p}\bar{q}}^n(T) \pl \le \pl \pq^n(T) \pl . \]
For $2 \le q < \infty$, $T \in \leb(X,Y)$ and $\nen$ the Rademacher
(gaussian)
Cotype $q$ norm with respect to $n$-vectors is defined by
\[ C_q^n(T) \pll (\p\gc_q^n(T)\p) \pl :=\pl \sup\{ \kla \summ_1^n \noo
Tx_k \rrm^q \mer^{\frac{1}{q}}
\mitt \kla \intt_{\Om} \noo \summ_1^n v_k x_k \rrm^2 \p
d\mu\mer^{\frac{1}{2}} \p\le \p 1\}\pl ,\]
where $(v_k)_1^n$ is a sequence of independent Bernoulli (gaussian)
variables on a probability space $(\Om,\mu)$.
An operator is said to be of Rademacher (gaussian) cotype $q$ if the
corresponding norm
\[ C_q(T) \pl :=\pl \sup_{\nen} C_q^n(T) \pl \kla \gc_q(T) \pl :=\pl
\sup_{\nen} \gc_q^n(T) \mer \]
is finite. For further information and the relation between gaussian
coype and
(q,2)-summing operators see for example \cite{TOJ}.

%%%%%%%%%%%%%%%%%%
%%Positive results
%%%%%%%%%%%%%%%%%%%

\section{Positive Results}

We start with the

{\bf Proof of theorem 2:} $\bf "\le" $ Let $x_1,..,x_n \in X$ with
\[ \sup_{\noo \al \rrm_{q'} \le 1} \noo \summ_1^n\al_k x_k \rrm \pl =
\pl
 \sup_{\noo x^*\rrm_{X^*}} \kla \summ_1^n \bet \langle x_k,x^*\rangle
 \rag^q \mer^{1/q}
\pl \le \pl 1 .\]
Therefore the operator $V \p:= \p \summ_1^n e_k \otimes x_k : \ell_{q'}
\nach X$ is of norm 1. By the equality case of H\"older's inequality we
obtain
\for
\kla \summ^n_1 \noo Tx_k \rrm^p \mer^{1/p} &=& \sup_{\noo \sigma\rrm_s
\le 1}
\kla \summ_1^n\kla \bet \sigma_k\rag \noo Tx_k \rrm \mer^r
\mer^{1/r}\\[+0.3cm]
&\le& \sup_{\noo \sigma\rrm_s \le 1} \pr^n(TV\ds) \pl .
\mel

$\bf"\ge"$ By the maximality of the norms $\pr^n$ there is no
restriction
to
assume $\ds:\lif^m \nach \ell_{q'}^m$ and $V\p:\p \ell_{q'}^m \nach X$
with
$\noo \sigma \rrm_s,\noo V \rrm \p \le \p 1$. Now let $U: \lif^n \nach
 \lif^m$ an operator of norm 1. By an observation of Maurey \cite{MAU}
the extreme points of such operators are of the form
\[ U \pl = \pl  \summ^n_1 e_k \otimes g_k \pl ,\]
with the $g_k$'s are of norm 1 and have disjoint support. Since we have
to
estimate the convex expression
\[ \kla \summ_1^n \noo TV\ds U(e_k)\rrm^r \mer^{1/r}\]
we can assume the that $U$ is of this form. We define $\tau$ and
$J:\ell_{q'}^n \nach \ell_{q'}^m$ by
\[ \tau_k \pl := \pl \noo \ds(g_k) \rrm_{q'} \pl \mbox{and} \pl
J\pl := \pl \summ_1^n e_k \otimes \frac{\ds(g_k)}{\noo \ds(g_k)
\rrm_{q'}}
\pl .\]
The operator J is of norm at most 1. Since $\ds$ is obviously
s1-summing
we have
\[ \noo \tau \rrm_s \pl \le \pl \pi_{s1}(\ds) \noo U \rrm
\pl \le \pl \noo \sigma \rrm_s  \pl \le \pl 1 \pl . \]
Therefore we obtain by H\"older's inequality
\for
\kla \summ_1^n \noo TV\ds U(e_k)\rrm^r \mer^{1/r}
 &=& \kla \summ_1^n \noo T V(\frac{\ds (g_k)}{\noo \ds (g_k)\rrm_{q'}}
 \tau_k)\rrm^r \mer^{1/r}\\
&=&\kla \summ_1^n \kla \noo TVJ(e_k)\rrm \p \bet \tau_k \rag \mer^r
\mer^{1/r}\\
&\le& \kla \summ_1^n \kla \noo TVJ(e_k) \rrm \mer^p \mer^{1/p} \pl \noo
\tau \rrm_s\\[+0.3cm]
&\le& \pq^n(T) \noo VJ \rrm \pl \noo \tau \rrm_s \pl \le \pl
\pq^n(T)\pl . \\[-1.5cm]
\mel \hfill $\Box$\hz

Now we will formulate a generalization of Jameson's lemma \cite{JAM}
which he proved in the case $q=1$, $p=2$.
\hz

\begin{lemma}[Jameson]\label{james} Let $1\le q <p \le \infty$ and
$T \in \leb(X,Y)$ a q-summing operator then
\[ \pi_{pq}(T) \pl \le \pl 2^{1/p} \pi_{pq}^n(T) \]
where
\[ n \pl \le \pl \kla 2^{1/p} \frac {\pi_q(T)}{\pq(T)}
\mer^{\frac{1}{1/q-1/p}} \pl .\]
\end{lemma}\hz

{\bf Proof:} Let us assume $\pi_q(T) \p=\p1$.
For $\eps > 0$ let $x_1,..,x_N$ in X with

\[ \sup_{\noo x^*\rrm_{X^*}\le1} \kla \summ_1^N \bet \langle
x_k,x^*\rangle \rag^q \mer^{1/q}
\pl \le \pl 1 \pll \mbox{and}
\pll (1-\eps) \pq^p(T) \pl \le \pl \summ_1^N \noo Tx_k \rrm^p \pl . \]
Furthermore, we assume $\noo Tx_k \rrm$ nonincreasing. For
$0 < \delta $
we choose $n \le N$ minimal such that $\noo Tx_k\rrm \le \delta$ holds
for all
$k> n$. Then we have $n \le \delta^{-q}$ because of
\[ n \delta^q \pl \le \pl \pi_q^q(T) \pl \le \pl 1 \pl . \]
If $\delta^{p-q}\p \le\p \frac{1}{2}(1-\eps) \pq^p(T)$ it follows that
\for
(1-\eps) \pq^p(T) &\le& \summ_1^N \noo Tx_k \rrm^p\\
&\le& \summ_1^n \noo Tx_k\rrm^p + \delta^{p-q}  \summ_{n+1}^N \noo
Tx_k\rrm^q\\
&\le& \summ_1^n \noo Tx_k\rrm^p + \delta^{p-q} \pi_q^q(T) \\
&\le&  \summ_1^n \noo Tx_k\rrm^p + \frac{1}{2}(1-\eps) \pq^p(T)\pl .
\mel
This means
\[ (1-\eps)^{1/p} \pq(T) \pl \le \pl 2^{1/p} \pq^n(T) \pl .\]
Letting $\eps$ to zero we find an $\nen$ with
\[ n \pl \le \pl \kla \frac{\pq^p(T)}{2} \mer^{\frac{-q}{p-q}}
\pl = \pl \kla 2^{1/p} \frac {\pi_q(T)}{\pq(T)}
\mer^{\frac{1}{1/q-1/p}} \pl .\hoch \]

\hfill $\Box$ \hz

\begin{samepage}
\begin{rem} \label{crem} Exactly the same argument shows that for every
operator $T \in \leb(X,Y)$
which is of (Rademacher) Cotype $2$ one has
\[ C_q(T) \pl \le \pl 2^{\frac{1}{q}} \pl C_q^n(T) \pl, \]
where $\nen$ satisfies
\[ n \pl \le \pl \kla 2^{1/q} \frac {C_2(T)}{C_q(T)}
\mer^{\frac{1}{1/q-1/2}} \pl .\]
\end{rem}\hz

{\bf Proof of theorem \ref{Rcot}:} Let E be a $n$-dimensional Banach
space. According to
Jameson's lemma we want to compare the cotype $2$ norm with the cotype
$q$ norm
via the gaussian cotype. It is well-known \cite[theorem 3.9.]{PS} that
the Rademacher cotype $2$
can be estimated by the gaussian cotype $2$ norm in the following way.
\[ C_2(Id_E) \pl \le \pl c_0 \pl \gc_2(Id_E) \pl \sqrt{1+\ln
\gc_2(Id_E)}\pl .\]
Using the inequalities $\gc_2(Id_E)\p\le\p\sqrt{2}\p
n^{\frac{1}{2}-\frac{1}{q}}\p \gc_q(Id_E)$
and $\gc_q(Id_E)\p\le\p n^{\frac{1}{q}}$, see \cite{TOJ}, we obtain
\for
 C_2(Id_E) &\le& c_0 \pl \gc_2(Id_E) \pl \sqrt{1+\ln \gc_2(Id_E)}\\
 &\le& \p c_0 \pl n^{\frac{1}{2}-\frac{1}{q}} \pl \gc_q(Id_E)\pl
 \sqrt{2+\ln n}
\mel
Combining this estimate with Jameson's lemma, more precisely the remark
above,
we see that there is a constant $c_1>0$ such that
\[ C_q(Id_E) \pl \le \pl 2^{\frac{1}{q}} \pl  C_q^m(Id_E)\pl ,\]
with
\for
m &\le& n \pl (c_1 \p \sqrt{1+\ln
n}\p)^{\frac{1}{1-\frac{2}{q}}}\\[-1.3cm]
\mel \hfill $\Box$ \end{samepage}\hz

In order to apply Jameson's lemma an appropriate estimate of the
1-summing norm by the r1-summing norm is needed.

\begin{lemma}\label{limit}
Let $1 \le r \le \infty$, $\nen$ and $T \in \leb(X,Y)$ an
operator of rank n. Then we have

\[ \pi_1(T) \pl \le \pl c_0 \pi_{r1}(T) \left \{
\begin{array}{l@{\quad}l}
    \kla \frac{1}{r}-\frac{1}{2}\mer^{-1/2} n^{1/2} & for \pl 1\le r<
    2\\[+0.2cm]
    \kla n(1+\ln n) \mer^{1/2} & r=2 \\[+0.2cm]
    \kla \frac{1}{2}-\frac{1}{r}\mer^{-1/r'} n^{1/r'} & for\pl
    2<r<\infty\pl .

 \end{array}\right. \]
\end{lemma}

{\bf Proof:} We may assume $T \in \leb(\lif,F)$ with dim$F$ = $n$. The
inequality $\pi_2(S: F \nach \lif) \p \le \p \sqrt{n}\noo S\rrm$, see
\cite{TOJ}, implies with Tomczak-Jagermann's inequality
\[ \pi_1(T) \pl \le \pl \iota_1(T) \pl \le \pl \sqrt{n}\p \pi_2(T) \pl
\le \pl
 \sqrt{2n}\p \pi_2^n(T) \pl .\]
For $2<r<\infty$ we deduce from Maurey's theorem, see \cite{TOJ}
\[ \pi_2^n(T) \pl \le \pl n^{1/2-1/r} \pi_{r2}^n(T)\pl \le \pl
  c_0 \kla \frac{1}{2}-\frac{1}{r}\mer^{-1/r'} n^{1/2-1/r} \pi_{r1}(T)
  \pl . \]
For $r=2$ we choose $2<\bar{r}<\infty$ with
$\frac{1}{2}-\frac{1}{\bar{r}}
 = \frac{1}{2+2\ln n}$. With $\pi_{\bar{r}1}\p\le\p\pi_{21}$ we obtain
\[ \pi_2^n(T) \pl \le \pl 2\pl e^2 \p c_0 (1+\ln n)^{1/2}\p
\pi_{21}^n(T) \pl . \]
In the case $1< r<2$ we use the other version of Maurey's
theorem, see again \cite{TOJ}, to deduce
\[\pi_2(T) \pl \le \pl c_0 \kla \frac{1}{r}-\frac{1}{2}\mer^{-1/2}
\pi_{r1}(T) \pl . \]
Combining the last three estimates with the first one gives the
assertion.
\hfill $\Box$ \hz

Now we can give the\hz

{\bf Proof of theorem 1:} First we prove the theorem in the case $q=1$,
hence $p'=r$.
From Jameson's lemma \ref{james} and lemma \ref{limit} we deduce for an
operator
$T$ of rank $n$
\[ \pi_{p1}(T) \pl \le \pl 2^{1/p} \pi_{p1}^m(T)\pl, \]
where
\[ m \pl \le \pl (2\p c_0)^r \left \{ \begin{array}{l@{\quad}l}
    \kla \frac{1}{2}-\frac{1}{p}\mer^{-1} n & for \pl 2<p<\infty\\
    n(1+\ln n) & for \pl r=2 \\
    \kla \frac{1}{p}-\frac{1}{2}\mer^{-r/2} n^{r/2} & for \pl 1\le p<
    2\pl.\\[+0.3cm]
 \end{array}\right. \]
An elementary computation shows that for all $\alpha,\p c \ge 1$ one
has
\[ \pi_{p1}^{[c \alpha]}(T) \pl \le \pl (4c)^{1/p}
\pi_{p1}^{[\alpha]}(T)
\pl . \]
Hence we get
\[\pi_{p1}(T) \pl \le \pl (16 \p c_0)^{r-1} \left \{\begin{array}{l
@{\quad} l}
 \kla \frac{1}{r}-\frac{1}{2} \mer^{-1/r'} \pl \pi_{p1}^n(T) & for \pl
 1 \le r <2\\[+0.2cm]
 \pi_{21}^{[n(1+\ln n)]} (T) & for \pl r=2\\[+0.2cm]
 \kla \frac{1}{2}-\frac{1}{r} \mer^{-(r-1)/2} \pl
 \pi_{p1}^{[n^{r/2}]}(T)
	    & for \pl 2<r<\infty \pl .
\end{array}\right.\]
For an arbitrary $1\le q \le p\le \infty$ we define $\bar{p}=r'$. Since
we
have $\frac{1}{\bar{p}}=\frac{1}{p}+\frac{1}{q'}$ we can deduce from
theorem
\ref{quot} and the inequalities above
\for
\pq(T) &=& \sup \{ \p \pi_{\bar{p}1}\kla TV\ds\mer \p \mid \p
    V:\ell_{q'} \nach X,\pl \ds:\lif \nach \ell_{q'}, \pl
    \noo \si \rrm_{q'}, \noo V \rrm\p \le \p 1 \}\\
&\le& c_r \sup \{ \p \pi_{\bar{p}1}^{m(r,n)} \kla TV\ds\mer \p \mid \p
 V:\ell_{q'} \nach X,\pl \ds:\lif \nach \ell_{q'}, \pl \noo \si
 \rrm_{q'},
 \noo V \rrm\p \le \p 1 \}\\
&=& c_r \pq^{m(r,n)}(T) \pl ,
\mel
where $m(r,n)=n$, $m(r,n)=[n(1+\ln n)]$, $m(r,n)=[n^{r/2}]$ for
$r<2$, $r=2$, $2<r$, respectively.\hspace*{\fill} $\Box$\hz

\begin{rem} The polynomial order of the vectors needed to compute
the pq-summing norm can be improved for several choices of p and q,
because they are close enough to the 2. Let $1\le q\le p\p, r\le\infty$
with
$\frac{1}{q}=\frac{1}{p}+\frac{1}{r}$. Then for all operators T of rank
n
one has
\[ \pq(T) \pl \le \pl (c_0)^{r/p} \left \{ \begin{array}{l@{\quad}l}
   \pq^{[n^{1+r(1/q-1/2)}]}(T) & for \pl 1\le q\le 2 \pl
		   \mbox{and}\pl 2\le r \le q'\\[+.2cm]
   \pq^{[n^{r(1/2-1/p)}]}(T) & for \pl 2 \le q \le \infty \pl.
\end{array} \right.\]
\end{rem}

{\bf Proof: First case:}
We choose $2\le s \le \infty$ such that
$\frac{1}{q}-\frac{1}{p}=\frac{1}{2}-\frac{1}{s}$.
By a result of Carl, \cite{CAR}, we have together with
Tomczak-Jaegermann's and Kwapien's inequality
\for
\pi_q(T)&\le& n^{1/q-1/2}\pl \pi_2(T) \
	    \pl \le \pl \sqrt{2} \pl n^{1/q-1/2} \pl
	    \pi_2^n(T)\\[+0.3cm]
&\le& \sqrt{2} \pl n^{1/q-1/2}\pl n^{1/2-1/s} \pl \pi_{s2}^n(T)
\pl \le \pl \sqrt{2} \pl n^{2/q-1/2-1/p} \pl \pq^n(T) \pl .
\mel
By Jameson's lemma \ref{james} and the elementry estimate in the prove
above
we obtain
\[ \pq(T) \pl \le \pl 8^{1/p} \pl
2^{\frac{r}{p}(\frac{1}{p}+\frac{1}{2})}
   \pl \pq^{[n^{1+r(1/q-1/2)}]} \pl . \]

{\bf Second case:} From Kwapien's and Tomczak-Jaegermann's inequality
we deduce
\for
\pi_q(T)&\le& \pi_2(T) \pl \le \pl \sqrt{2}\pl  \pi_2^n(T)\\[+0.3cm]
 &\le& \sqrt{2} \pl n^{1/2-1/p} \pl \pi_{p2}^n(T) \pl \le \pl
	\sqrt{2}\pl n^{1/2-1/p}\pl \pq(T) \pl .
\mel
Again with Jameson's lemma \ref{james} this implies the assertion.
\hfill $\Box$\hz

At the end of this chapter we would like to note the following

\begin{cor} Let $1 \le r \le \infty$, K a compact Hausdorff space and
 $T \in \leb(C(K),Y)$ of rank n. Then we have
\[ \pi_p(T) \pl \le \pl c_p (1+ \ln n)^{1/p'} \left \{
\begin{array}{l @{\quad}l}
\pi_p^n(T) & for \pl 2<p<\infty\\[+0.2cm]
\pi_p^{[n^{p'/2}]}(T) & for \pl 1<p<2 \pl .
\end{array} \right. \]
\end{cor}

{\bf Proof:} Using a result of Carl and Defant, see \cite{CAD}, and
theorem 2 we deduce
\for
\pi_p(T) &\le& c_0 \pl (1+ \ln n)^{1/p'} \pl \pi_{p1} (T) \\
& \le & c_p \pl (1+ \ln n)^{1/p'} \pl \pi_{p1}^{[n^{\max(1,p'/2}]}(T)
\pl \le \pl  c_p \pl (1+ \ln n)^{1/p'} \pl \pi_p^{[n^{\max(1,p'/2}]}(T)
\pl .\\[-1.2cm]
\mel \hfill $\Box$ \hz

%%%%%%%%%%%
%%Examples
%%%%%%%%%%
\section{Examples}
\setcounter{lemma}{0}

By the positve results of the previous section a polynomial groth can
only appear if $\frac{1}{q}-\frac{1}{p}>\frac{1}{2}$. Therefore we
define
for $1 \le q \le 2$ the critical value $p_q$ by
$\frac{1}{p_q}=\frac{1}{q}-\frac{1}{2}$. In
the sequel limit orders of pq-summing operators are of particular
interest. We intensively use the results of Carl, Maurey and Puhl, see
\cite{CMP}. The next lemma is implicitely contained there but we
reproduce the easy interpolation argument.

\begin{lemma}\label{int} Let $0<\tet<1$, $1 \le q,\p r \le 2$ and $q\le
p$
with $\frac{1}{r}=\frac{1-\tet}{q}+\frac{\tet}{2}$ and
$\frac{1}{p}=\frac{1}{q}-\frac{\tet}{2}$. Then one has
\[ \pq(\iota:\ell_{q'}^n\nach \ell_r^n) \pl \le \pl n^{1/p}\pl . \]
\end{lemma}

{\bf Proof:} Clearly one has, see \cite{PIE},
\[ \pi_q(\iota : \ell_{q'}^n \nach \ell_q^n) \pl \le \pl
   \pi_q(\iota : \lif^n \nach \ell_q^n) \pl \le \pl n^{1/q} \pl . \]
With Kwapien's inequality, $\pi_{qp_q} \p \le \p \pi_{21}$, we deduce
from the Orlicz property of $\ell_2$
\[ \pi_{qp_q}(\iota : \ell_{q'}^n \nach \ell_2^n) \pl \le \pl
 \pi_{21}(Id_{\ell_2^n})
 \pl \noo \iota : \ell_{q'}^n \nach \ell_2^n \rrm \pl \le \pl
 n^{1/q-1/2} \pl . \]
By interpolation, namely $[\ell_q(\ell_q),\ell_{p_q}(\ell_2)]_{\tet}
=\ell_p(\ell_r)$ see \cite{BEL}, this means
\for
\pq(\iota:\ell_{q'}^n\nach \ell_r^n) &\le&
 n^{\frac{1-\tet}{q}} \pl n^{\tet(\frac{1}{q}-\frac{1}{2})}
\pl = \pl n^{1/p} \pl . \\[-1.5cm]
\mel
\hfill $\Box$\hz

%%%%%%%%%%%%%%%%%%%%%%%%%%%%%%%%%%%%%%%%%%%
Now we can construct the counterexamples

\begin{prop} \label{tetversion} Let $0<\tet<1$, $1 \le q \le 2\le s \le
\infty$ and $q\le t$
with $\frac{1}{s}=\frac{1-\tet}{q'}+\frac{\tet}{2}$ and
$\frac{1}{t}=\frac{1}{q}-\frac{\tet}{2}$. For all $\nen$ there exists
an operator $T \in \leb(\ell_{q'},\ell_2^n)$ wich satisfies
\[ \frac{\pq^k(T)}{\pq(T)} \pl \le \pl c_0 \sqrt{s}
\kla \frac{k}{n^{s/2}} \mer^{1/p-1/t} \]
for all $q \le p \le t$ and $\ken$.
\end{prop}

{\bf Proof:} Let $m=[n^{s/2}]$ and $A:\ell_{q'}^m \nach  \ell_2^n$
be a random matrix with entries $\pm 1$,
a so called Benett
matrix. Obviously we have
\[ \pq(A) \pl \ge \pl m^{1/p}\pl n^{1/2} \pl \ge \pl \frac{1}{2} \pl
n^{1/2+s/2p} \pll
\mbox{for all}\pl  q \le p \pl . \]
We will see that this estimate is sharp for some indices p. By
\cite[Lemma 5]{CMP} one has
\[ \noo A: \ell_{s'}^m \nach  \ell_2^n \rrm \pl \le \pl
 c_0 \sqrt{s} \max\{n^{1/2},m^{1/s}\} \pl \le \pl c_0 \sqrt{s\p n}\pl .
 \]
Since $\frac{1}{s'}=\frac{1-\tet}{q}+\frac{\tet}{2}$ we deduce
from Lemma \ref{int}
\for
\pi_{tq}(A) &\le& \pi_{tq}(\iota:\ell_{q'}^m \nach \ell_{s'}^m)
  \pl \noo A: \ell_{s'}^m \nach  \ell_2^n \rrm      \\
 &\le& \pl c_0 \p \sqrt{s} \pl n^{1/2}\pl m^{1/t}
\pl \le \pl  c_0 \p \sqrt{s} \pl n^{1/2+s/2t} \pl .
\mel
Therefore we obtain for arbitrary $q \le p \le t$, $\ken$
\for
\pq^k(A) &\le& k^{1/p-1/t} \pl \pi_{tq}^k(A)\\
&\le& c_0 \p \sqrt{s} \pl\pl k^{1/p-1/t} \pl\pl n^{1/2+s/2t}\\
&\le& c_0 \p \sqrt{s} \pl \kla \frac{k}{n^{s/2}} \mer^{1/p-1/t}
 \pl \pq(A) \pl . \\[-1.3cm]
\mel \hfill $\Box$ \hz
%%%%%%%%%%%%%%%%%%%%%%%%%

\begin{cor}\label{order} For $1 \le q < 2$ and $q\le  p < p_q$ let
$2 < s\le \infty$ defined by
$\frac{1}{s}=\frac{1}{q'}+(1-\frac{2}{q'})(\frac{1}{q}-\frac{1}{p})$.
For any sequence $k_n$, $c>0$ satisfying
\[ \pq(T) \pl \le \pl C \pq^{k_n}(T)  \pll \pl \mbox{for all T of rank
n}\]
there is a constant $c_1$ with
\[ n^{s/2} \pl \le \pl c_1 e^{c_1 \p \sqrt{1+\ln n}} \pl k_n  \pl . \]
\end{cor}

{\bf Proof:} We define $\ttett:=2(\frac{1}{q}-\frac{1}{p}) < 1$.
For $\eps < 1-\ttett$ we set $\tet:=\ttett+\eps$ and choose $2 \le v
\le s$,
$p \le t \le p_q$ with
$\frac{1}{v}=\frac{1-\tet}{q'}+\frac{\tet}{2}$,
$\frac{1}{t}=\frac{1}{q}-\frac{\tet}{2}$. Now let us consider
the quotient $d_n\p := \p n^{s/2}k_n^{-1}$. From Proposition
\ref{tetversion} we deduce with an elementary computation
\for
\frac{1}{C} &\le& c_0\p\sqrt{s} \pl\kla
\frac{k_n}{n^{v/2}}\mer^{1/p-1/t}
 \pl \le \pl c_0 \p\sqrt{s} \pl d_n^{-\eps/2} \pl n^{(s-v) \eps/4}  \\
&\le& c_0 \p\sqrt{s}\pl d_n^{-\eps/2}\pl n^{\eps^2 s^2(1/8-1/4q')} \pl
.
\mel
Insetting $\eps =\frac{1-\ttett}{2}$ yields a constant $c_2$ such that
\[ \ln d_n \pl \le \pl c_2 + \kla s^2 (1-\ttett) \p
 (\frac{1}{8}-\frac{1}{4q'}) \mer \ln n \pl . \]
Therefore there exists an $n_0 \in \nz$ such that for all $n \ge n_0$
we can choose $\eps_n := \frac{\ln d_n}{s^2(1/2-1/q')\ln n} \p <
1-\ttett$.
An elementary computation gives
\[ \frac{1}{C}\pl \le \pl c_q \exp\kla-
\frac{(\ln d_n)^2}{2s^2(1/2-1/q')(\ln n)}\mer \pl . \]
This is only possible if there exists a constant $c_3$ depending on s,
q and
C such that
\for
 \frac{n^{s/2}}{k_n} \pl = \pl d_n &\le&\exp\kla c_3 \sqrt{1+ \ln n}
 \mer \pl .\\[-1.3cm]
\mel \hfill $\Box$
%%%%%%%%%%%%%%%%%%%

\begin{rem} For $q=1$ the above results can be slightly improved. For
$1 \le p < 2 \le s < p'\le \infty$ there exists an operator
$T \in \leb(\lif^{[n^{s/2}]},\ell_2^n)$ such that for all $\ken$
\[ \pi_p^k(T)\pl \le \pl c_0 \sqrt{s} \p \kla \frac{k}{n^{s/2}}
\mer^{1/s-1/p'} \pl \pi_{p1}(T) \pl . \]
In particular, the inequality
\[ \pi_{p1}(T) \pl \le \pl C\p (1+\ln n)\pl  \pi_p^{k_n}(T) \pla
\mbox{for all
T of rank n}\]
can only be satisfied, if
\[ n^{p'/2} \pl \le \pl \bar{c}\pl e^{\bar{c}\p\sqrt{1 +\ln n}}\pl k_n
\pl . \]
This answers a conjecture of Carl and Defant. They suggested
\[ \pi_p(T) \pl \le \pl c_p (1+\ln n)^{1/p'} \pi_{p1}^n(T) \]
for all operators $T\in \leb(C(K),Y)$ of rank n, which turns out to be
false.
Furthermore, we recover the exponential order of vectors for $\pi_1$.
More precisely, for all $n,k \in \nz$ there is an operator
$T \in \leb(\lif^{[n^{1+\ln k}]},\ell_2^n)$ with
\[ \pi_1^k(T) \pl \le \pl c_0\pl \sqrt{\frac{1+\ln k}{n}} \pl \pi_1(T)
\pl . \]
\end{rem}\hz

{\bf Proof:} Inspecting the proof of proposition \ref{tetversion} we
take
a Benett matrix $A: \lif^{[n^{s/2}]} \nach  \ell_2^n$, whose
p1-summing norm satisfies
\[ \pi_{p1}(T) \pl \ge\pl \frac{1}{2}\pl n^{1/2+s/2p} \pl . \]
On the other hand
\for
 \pi_p^k(A)&\le& k^{1/p-1/s'}\pl \pi_{s'}(A)\\
&\le& k^{1/p-1/s'}\pl  \pi_{s'}(\iota : \lif^{[n^{s/2}]} \nach
\ell_{s'}^{[n^{s/2}]})
\pl \noo A: \ell_{s'}^{[n^{s/2}]} \nach \ell_2^n \rrm \\
&\le& c_0\p \sqrt{s}\pl k^{1/p-1/s'}\pl n^{s/2s'}\pl  n^{1/2} \\
&\le& 2 c_0 \p \sqrt{s} \pl \kla \frac{k}{n^{s/2}}
\mer^{1/p-1/s'} \pl \pi_{p1}(T) \pl .
\mel
The logarithmic factor does not affect the calculation in the
proof of corollary \ref{order}. For the last assertion we note that
$p'=\infty$ and therefore the choice $s = 2(1+\ln k)$ implies the
assertion.

\hfill $\Box$ \hz

%%%%%%%%%%%%%%

Now we will give the\hz

{\bf Proof of Theorem 4:}  In the case $1 \le q<2$ this follows
immiadetly from corollary \ref{order}. We only have to note that
for all $\eps >0$ there is a constant $C_{\eps}$ with
\[ e^{\bar{c}\sqrt{1+\ln n}}\pl \le \pl C_{\eps} n^{\eps}\pl . \]
Now let $2 < q < \infty$. With the help of Benett matrices it was
shown in \cite{CMP} that for $2<q<\infty$
\[ n^{q/2p} \pl \le \pl c_0 \p \sqrt{q} \pl \pq (id_{\ell_2^n})\pl . \]
Hence we get
\for
\pq^k(id_{\ell_2^n}) &\le& k^{1/p} \noo id_{\ell_2^n}\rrm \pl \le \pl
k^{1/p}\\
&\le& c_0 \p\sqrt{q} \pl \kla \frac{k}{n^{q/2}} \mer^{1/p} \pl
\pq(id_{\ell_2^n})\pl .\\
\mel
Therefore every sequence $k_n$ with $\pq(T) \le C \p\pq^{k_n}(T)$
must satisfy
\for
 n^{q/2} &\le& \kla C\p c_0 \p\sqrt{q} \mer^p \pl k_n \pl .\\[-1.4cm]
\mel
\hfill $\Box$
%%%%%%%%%%%

\begin{rem} For operators defined on $n$-dimensional Banach spaces
the results of {\rm\cite{JOS}} and {\rm\cite{DJ}} imply that the
pq-summing norm
can be
calculated with $n^{q/2} \p (1+\ln n)$ many vectors. Therefore
the order in the proof of the proposition above is quite correct.
\end{rem}
%%%%%%%%%%%%
%%%%%%%%%%%%References
%%%%%%%%%%%%%

\hz
1991 Mathematics Subject Classification: 47B10, 47A30,46B07.
\begin{quote}
Martin Defant

Marius Junge

Mathematisches Seminar der Universit\H{a}t Kiel

Ludewig-Meyn-Str. 4

24098 Kiel

Germany

Email: nms06@rz.uni-kiel.d400.de
\end{quote}

\end{document}